\documentclass[12pt,oneside,reqno]{amsart}
\usepackage{mathrsfs}
\usepackage{graphics}
\usepackage{float}
\usepackage{subfigure}
\usepackage{color, amssymb}
\newcommand{\dif}{\mathrm{d}}
\newcommand{\me}{\mathrm{e}}
\newcommand{\be}{\begin{eqnarray}}
\newcommand{\ee}{\end{eqnarray}}
\newcommand{\ce}{\begin{eqnarray*}}
\newcommand{\de}{\end{eqnarray*}}
\newtheorem{theorem}{Theorem}[section]
\newtheorem{lemma}[theorem]{Lemma}
\newtheorem{remark}[theorem]{Remark}
\newtheorem{definition}[theorem]{Definition}
\newtheorem{proposition}[theorem]{Proposition}
\newtheorem{example}[theorem]{Example}
\newtheorem{corollary}[theorem]{Corollary}

\newcommand{\cL}{{\mathcal{L}}}

\newcommand{\PX}{{\Bbb{P}}}

\def\e{\varepsilon}

\def\a{\alpha}

\def\[{{\Big[}}
\def\]{{\Big]}}
\def\<{{\langle}}
\def\>{{\rangle}}
\def\({{\Big(}}
\def\){{\Big)}}

\def\no{\nonumber}
\def\bt{\begin{theorem}}
\def\et{\end{theorem}}
\def\bl{\begin{lemma}}
\def\el{\end{lemma}}
\def\br{\begin{remark}}
\def\er{\end{remark}}
\def\bx{\begin{Examples}}
\def\ex{\end{Examples}}
\def\bd{\begin{definition}}
\def\ed{\end{definition}}
\def\bp{\begin{proposition}}
\def\ep{\end{proposition}}
\def\bc{\begin{corollary}}
\def\ec{\end{corollary}}

\def\cC{{\mathcal C}}

\def\cF{{\mathcal F}}

\def\cL{{\mathcal L}}

\def\mE{{\mathbb E}}

\def\mR{{\mathbb R}}

\def\geq{\geqslant}
\def\leq{\leqslant}

\begin{document}

\allowdisplaybreaks

\title{Asymptotic methods for stochastic dynamical systems with small non-Gaussian L\'evy noise*}

\author{Huijie Qiao$^1$ and Jinqiao Duan$^2$}



\thanks{*This work was partially supported by the NSFC (No. 11001051, 10971225 and 11028102) and the NSF grant DMS-1025422.}


\date{\today}

\dedicatory{1. Department of Mathematics,
Southeast University\\
Nanjing, Jiangsu 211189,  China\\
hjqiaogean@seu.edu.cn\\
2. Institute for Pure and Applied Mathematics, University of California\\
Los Angeles, CA 90095\\ \& \\Department of Applied Mathematics, Illinois Institute of Technology\\
Chicago, IL 60616, USA\\
duan@iit.edu}

\begin{abstract}
The goal of the paper is to analytically examine escape probabilities for dynamical systems
driven by symmetric $\alpha$-stable L\'evy motions. Since escape probabilities are solutions of a
type of integro-differential equations (i.e., differential equations with nonlocal interactions),
asymptotic methods are offered to solve these equations to  obtain escape probabilities when noises
are sufficiently small. Three examples are presented to illustrate the asymptotic methods, and
asymptotic escape probability is compared with numerical simulations.
\end{abstract}

\maketitle \rm

\section{Introduction}
Stochastic dynamical systems are mathematical models from complex phenomena in biological, geophysical,
physical and chemical sciences, under random fluctuations.
Unlike the situation for deterministic dynamical systems, an orbit  of a stochastic system
could vary wildly from one sample to another. It is thus desirable to have efficient tools
to quantify stochastic dynamical behaviors. The escape probability is such a tool.

Non-Gaussian random fluctuations are widely observed in various areas such as physics, biology, seismology,
electrical engineering and finance \cite{Woy, Koren, Mat}. L\'evy motions are a class of non-Gaussian
  processes whose sample paths are discontinuous in time.  For a dynamical system  driven by L\'evy motions,
almost all   orbits  are discontinuous in time. In fact, these orbits are c\`adl\`ag (right-continuous with
left limit at each time instant), i.e., each of these orbits has countable jumps in time. Due to these jumps, an orbit
could escape an open domain without passing through its boundary. In this case, the \emph{escape probability} is the
likelihood that an orbit, starting inside an open domain $D$, exits this domain first by landing in a target domain
$U$ in $D^c$ (the complement of domain $D$).

For brevity, in this paper we only consider scalar stochastic dynamical systems. Let $\{X_t, t\geq0\}$ be a real-valued Markov process defined on a complete filtered probability space
$(\Omega,\cF,\{\cF_t\}_{t\geq 0}, \PX)$. Let $D$ be an open domain in $\mR$. Define the \emph{exit time}
  \ce
 \tau_{D^c}:=\inf\{t>0: X_t\in D^c\},
 \de
where $D^c$ is the complement of $D$ in $\mR$. Namely, $\tau_{D^c}$ is the  first time when $X_t$ hits $D^c$.

When $X_t$ has c\`adl\`ag paths which have countable jumps in time, starting
at $x\in D$, the first hitting of $D^c$ may occur either on the boundary
$\partial D$ or somewhere in $D^c$. For this reason, we take a
subset $U$ of the closed set $D^c$, and define the likelihood that   $X_t$ exits
firstly from $D$ by landing in the target set $U$ as the escape probability from $D$ to
$U$, denoted by $p(x)$. That is,
\ce
p(x)=\PX \{X_{\tau_{D^c}}(x)\in U\}.
\de

If $X_t$ is a solution process of a dynamical system driven by a symmetric $\alpha$-stable L\'evy motion, by \cite{Liao, Kan},
the escape probability $p(x)$ solves the following Balayage-Dirichlet ``exterior" value problem:
\be
\left\{\begin{array}{l}
Ap=0,\;\; x\in D, \\
p|_{D^c}=\varphi,
\end{array}
\right.
\label{dirpro1}
\ee
where $A$ is the infinitesimal generator of $X_t$ and $\varphi$ is defined as follows
\ce
\varphi(x)=\left\{\begin{array}{l}
1, \quad x\in U,\\
0, \quad x\in D^c\setminus U.
\end{array}
\right.
\de

However, Eq.(\ref{dirpro1}) is usually an integro-differential equation and it is hard to
obtain exact representations for its solutions. Here we use asymptotic methods to examine
its solutions. More precisely, \textbf{(i)} for a dynamical system driven by a Brownian
motion combined with a symmetric $\alpha$-stable L\'evy motion, an asymptotic
solution of Eq.(\ref{dirpro1}), or escape probability $p(x)$ from $D$ to $U$, is obtained
by a regular perturbation method; \textbf{(ii)} for a dynamical system driven by a symmetric
$\alpha$-stable L\'evy motion alone, the escape probability $p(x)$ is obtained by a singular perturbation method.

This paper is arranged as follows. In Section \ref{prelim}, we introduce symmetric $\alpha$-stable
L\'evy motions and their infinitesimal generators. In Section \ref{sde-contjump}, a regular perturbation
method is  applied to examine escape probability for dynamical systems driven jointly by Brownian motion and symmetric $\alpha$-stable L\'evy motions. Escape probabilities for dynamical systems driven by symmetric $\alpha$-stable L\'evy motions alone
are studied in Section \ref{sde-jump} by a singular perturbation method.   Three examples are presented in
Section \ref{example}.

\section{Preliminaries}\label{prelim}

In this section, we recall  basic concepts and results that will
be needed throughout the paper.

\bd\label{levy}
A process $L_t=(L_t)_{t\geq0}$ with $L_0=0$ a.s. is a   L\'evy process or L\'evy motion if

(i) $L_t$ has independent increments; that is, $L_t-L_s$ is
independent of $L_v-L_u$ if $(u,v)\cap(s,t)=\emptyset$;

(ii) $L_t$ has stationary increments; that is, $L_t-L_s$ has the same
distribution as $L_v-L_u$ if $t-s=v-u>0$;

(iii) $L_t$ is stochastically continuous;

(iv) $L_t$ is right continuous with left limit.
\ed
The characteristic function for $L_t$ is
\ce
\mE\left(\exp\{izL_t\}\right)=\exp\{t\Psi(z)\}, \quad
z\in\mR.
\de
We only consider scalar L\'evy motions here. The function $\Psi: \mR\rightarrow\mathcal {C}$ is
called the characteristic exponent of the L\'evy process $L_t$. By
L\'evy-Khintchine formula, there exist a nonnegative number $Q$,
a measure $\nu$ on $\mR$ satisfying
\ce
\nu(\{0\})=0 ~\mbox{and}~ \int_{\mR\setminus\{0\}}(|u|^2\wedge1)\nu(\dif
u)<\infty,
\de
and also a real number $\gamma$ such that
\be
\Psi(z)&=&i\gamma z-\frac{1}{2}Qz^2+\int_{\mR\setminus\{0\}}
\big(e^{izu}-1-izu1_{|u|\leq1}\big)\nu(\dif
u).
\label{lkf}
\ee
The measure $\nu$ is called the L\'evy measure, $Q$ is the diffusion,
and $\gamma$ is the drift.

\medskip

We now introduce a special class of L\'evy motions, i.e., the symmetric $\alpha$-stable
L\'evy motions $L_t^\alpha$.

\bd\label{rid}
For $\alpha\in(0,2)$. A scalar symmetric $\alpha$-stable L\'evy motion $L_t^\alpha$ is a
L\'evy process with characteristic exponent
\ce
\Psi(z)=-|z|^\alpha, \quad z\in\mR.
\de
\ed

Thus, for a scalar symmetric $\alpha$-stable L\'evy motion $L_t^\alpha$, the diffusion $Q=0$,
the drift $\gamma=0$, and the L\'evy measure $\nu$ is given by
$$
\nu(\dif u)=\frac{C_{1,\alpha}}{|u|^{1+\alpha}}\dif u, \quad
C_{1,\alpha}=\frac{\a\Gamma((1+\a)/2)}{2^{1-\a}\pi^{1/2}\Gamma(1-\a/2)}.
$$

Let $\cC_0(\mR)$ be the space of continuous functions $f$ on $\mR$ satisfying
$\lim\limits_{|x|\rightarrow\infty}f(x)=0$ with norm $\|f\|_{\cC_0(\mR)}=\sup\limits_{x\in\mR}|f(x)|$.
Let $\cC^2_0(\mR)$ be the set of $f\in\cC_0(\mR)$ such that $f$ is twice
  differentiable and the first and second order   derivatives of $f$ belong
to $\cC_0(\mR)$. Let $\cL_\alpha$ be the infinitesimal generator of $L_t^\alpha$.
By \cite[Theorem 31.5]{s},
\ce
(\cL_\alpha f)(x)=\left\{\begin{array}{l}
\int_{\mR\setminus\{0\}}\big(f(x+u)-f(x)\big)\nu(\dif u), \quad\qquad\qquad\qquad 0<\alpha<1,\\
\int_{\mR\setminus\{0\}}\big(f(x+u)-f(x)
-f'(x)u1_{|u|\leq1}\big)\nu(\dif u), \qquad\alpha=1,\\
\int_{\mR\setminus\{0\}}\big(f(x+u)-f(x)
-f'(x)u\big)\nu(\dif u), \qquad\quad 1<\alpha<2,
\end{array}
\right.
\de
where $f\in\cC_0^2(\mR)$. For any $\e>0$, $\e L_t^\alpha$ is also a scalar
symmetric $\alpha$-stable L\'evy motion, and its L\'evy measure $\nu^\e(B)=\nu(\frac{1}{\e}B)$
for $B\in\mathscr{B}(\mR)$ (Borel $\sigma$-algebra on $\mR$). Thus, its infinitesimal generator
is
\ce
(\cL_\alpha f)(x)=\left\{\begin{array}{l}
\int_{\mR\setminus\{0\}}\big(f(x+u)-f(x)\big)\nu^\e(\dif u), \quad\qquad\qquad\qquad 0<\alpha<1,\\
\int_{\mR\setminus\{0\}}\big(f(x+u)-f(x)
-f'(x)u1_{|u|\leq1}\big)\nu^\e(\dif u), \qquad\alpha=1,\\
\int_{\mR\setminus\{0\}}\big(f(x+u)-f(x)
-f'(x)u\big)\nu^\e(\dif u), \qquad\quad 1<\alpha<2,
\end{array}
\right.
\de
Applying the representation of $\nu$, one can obtain that
\ce
(\cL_\alpha f)(x)=\left\{\begin{array}{l}
\e^\alpha\int_{\mR\setminus\{0\}}\big(f(x+u)-f(x)\big)\nu(\dif u), \quad\qquad\qquad\qquad 0<\alpha<1,\\
\e^\alpha\int_{\mR\setminus\{0\}}\big(f(x+u)-f(x)
-f'(x)u1_{|u|\leq1}\big)\nu(\dif u), \qquad\alpha=1,\\
\e^\alpha\int_{\mR\setminus\{0\}}\big(f(x+u)-f(x)
-f'(x)u\big)\nu(\dif u), \qquad\quad 1<\alpha<2,
\end{array}
\right.
\de

\section{Escape probability of a \textsc{SDE} with a Brownian motion and a symmetric $\alpha$-stable L\'evy motion}
\label{sde-contjump}

Let $\{W(t)\}_{t\geq 0}$ be a scalar standard $\cF_t$-adapted Brownian motion,
and $L_t^{\alpha}$ a scalar symmetric $\alpha$-stable L\'evy motion with $\alpha \in (0, 2)$ and independent
of $W_t$. Consider the following scalar stochastic differential equation, with the drift coefficient $b$, the
diffusion coefficient $\sigma$  and intensity $\e (>0)$ of L\'evy noise,
\be\left\{\begin{array}{l}
\dif X_{t}=b(X_{t})\,\dif t+\sigma(X_t)\,\dif W_t+\e\dif L_t^{\alpha},\\
X_0=x.
\end{array}
\right.
\label{sdebrolevy}
\ee
Assume that the  drift  $b$ and the  diffusion $\sigma(\neq0)$ satisfy the following conditions:

{\bf (H$_b$)}
there exists a constant $C_b>0$ such that for $x,y\in\mR$
$$
|b(x)-b(y)|\leq C_b|x-y|\cdot
\log(|x-y|^{-1}+\me);
$$

{\bf (H$_\sigma$)}
there exists a constant $C_\sigma>0$ such that for $x,y\in\mR$
$$
|\sigma(x)-\sigma(y)|^2\leq C_\sigma|x-y|^{2}\cdot
\log(|x-y|^{-1}+\me).
$$

Under {\bf (H$_b$)} and {\bf (H$_\sigma$)}, it is well known that there exists a
unique strong solution to Eq.(\ref{sdebrolevy})(see \cite{hqxz}). This solution will be
denoted by $X_t(x)$. By Theorem 3.3 in \cite{Kan}, the escape probability $p(x)$ for
$X_{t}(x)$ with $0<\alpha<1$, from $D=(A,B)$ to $U=[B,\infty)$, satisfies the following integro-differential equation
\be
b(x)p'(x)+\frac{1}{2}\sigma^2(x)p^{''}(x)+\e^\alpha\int_{\mR\setminus\{0\}}\left(p(x+u)-p(x)\right)\nu(\dif u)=0, \quad
x\in(A,B),
\label{intdif1}
\ee
and the `exterior' conditions
\begin{equation}\label{exterior}
p(x)|_{(-\infty,A]}=0, \qquad p(x)|_{[B,\infty)}=1.
\end{equation}

\medskip

We consider the  solution for $p(x)$, when $\e>0$ is sufficiently small.
Assume that $p(x)$ has the following regular expansion
\be
p(x)=p_0(x)+\e^\alpha p_1(x)+\e^{2\alpha}p_2(x)+\cdots.
\label{asym1}
\ee
Substituting (\ref{asym1}) into (\ref{intdif1}) and equating like powers of $\e$, we obtain a system of equations for the recursive determination of $p_j(x)$. The leading-order
equation for $p_0(x)$ is
\be
b(x)p_0^{'}(x)+\frac{1}{2}\sigma^2(x)p_0^{''}(x)=0, \qquad x\in(A,B)
\label{ode1}
\ee
with   boundary conditions
\be
p_0(A)=0, \qquad p_0(B)=1.
\label{orb1}
\ee
Using the boundary conditions, we solve Problem (\ref{ode1}) and (\ref{orb1}) to get
\ce
p_0(x)=\frac{\int_A^xe^{-\int_A^s\phi(u)\dif u}\dif s}{\int_A^Be^{-\int_A^s\phi(u)\dif u}\dif s},
\de
where $\phi(u):=2b(u)/\sigma^2(u)$.

Next, the equation for $p_1(x)$ is
\be
b(x)p_1^{'}(x)+\frac{1}{2}\sigma^2(x)p_1^{''}(x)+\int_{\mR\setminus\{0\}}\left(p_0(x+u)-p_0(x)
\right)\nu(\dif u)=0, \quad
x\in(A,B),
\label{ode2}
\ee
with boundary conditions
\be
p_1(A)=0, \qquad p_1(B)=1.
\label{orb2}
\ee
Set
\ce
g(x):=\int_{\mR\setminus\{0\}}\left(p_0(x+u)-p_0(x)\right)\nu(\dif u).
\de
Then Eq.(\ref{ode2}) is transformed into the following equation
\be
b(x)p_1^{'}(x)+\frac{1}{2}\sigma^2(x)p_1^{''}(x)+g(x)=0, \qquad
x\in(A,B).
\label{ode3}
\ee
By solving Problem (\ref{ode3}) and (\ref{orb2}) we get
\ce
p_1(x)&=&\int_A^xe^{-\int_A^s\phi(u)\dif u}\left(\int_A^s\frac{-2g(u)}{\sigma^2(u)}\cdot
e^{\int_A^u\phi(v)\dif v}\dif u\right)\dif s\\
&&-p_0(x)\int_A^Be^{-\int_A^s\phi(u)\dif u}\left(\int_A^s\frac{-2g(u)}{\sigma^2(u)}\cdot
e^{\int_A^u\phi(v)\dif v}\dif u\right)\dif s+p_0(x).
\de

Thus we have an asymptotic expression for escape probability, i.e., solution of Eq.(\ref{intdif1}), for $\e$ sufficiently small,
\begin{equation}
p(x)\approx p_0(x)+\e^\alpha p_1(x).
\end{equation}

By the same deduction as above, we could obtain an asymptotic expression for escape probability of
$X_{t}(x)$ with $1\leq\alpha<2$.

\section{Escape probability of a \textsc{SDE} with
a symmetric $\alpha$-stable L\'evy motion}\label{sde-jump}

Consider the following stochastic differential equation with  a symmetric
$\alpha$-stable L\'evy motion, with $1<\a<2$, on $\mR$
\be\left\{\begin{array}{l}
\dif X_{t}=b(X_t)\,\dif t+\e\dif L_t^{\alpha},\\
X_0=x,
\end{array}
\right.
\label{sdelevy}
\ee
where the drift $b$ satisfies {\bf (H$_b$)}.

By \cite[Theorem 3.1]{ttw},
Eq.(\ref{sdelevy}) has a unique solution   $X_{t}(x)$. From Theorem 3.3 in
\cite{Kan}, the escape probability $p(x)$, for $X_{t}(x)$ from $D=(A,B)$ to $U=[B,\infty)$, satisfies
the following integro-differential equation
\be
b(x)p'(x)+\e^\alpha\int_{\mR\setminus\{0\}}\big(p(x+u)-p(x)
-p'(x)u\big)\nu(\dif u)=0,\no\\
x\in(A,B),
\label{intdif2}
\ee
 with the `exterior' conditions
\be
&&p(x)|_{(-\infty,A]}=0,\label{boucon1}\\
&&p(x)|_{[B,\infty)}=1.
\label{boucon2}
\ee

We now try to construct an asymptotic solution of (\ref{intdif2}), (\ref{boucon1}) and (\ref{boucon2}) for sufficiently small $\e>0$. We
consider the following four different cases, depending on the dynamical behavior of the corresponding deterministic dynamical system $\dot{x}=b(x)$.

\medskip


{\bf Case 1:}   $b(x)>0$ for $x\in(A,B)$.  In this case the deterministic dynamical system $\dot{x}=b(x)$ has no equilibrium states and all orbits move to the right.

Thus it is reasonable to require that $p(x)\rightarrow1$ as $\e\rightarrow0$.
We assume that $p(x)$ has the following   expansion
\be
p(x)=p_0(x)+\e^\alpha p_1(x)+\e^{2\alpha}p_2(x)+\cdots.
\label{asym2}
\ee
Substituting (\ref{asym2}) into (\ref{intdif2}) and equating like powers of $\e$, we
obtain a system of equations for the recursive determination of $p_j(x)$. The leading-order
equation for $p_0(x)$ is
\ce
b(x)p_0^{'}(x)=0, \qquad x\in(A,B),
\de
and thus $p_0(x)=1$ for $x\in(A,B)$, because $p(x)\rightarrow1$ as $\e\rightarrow0$. Since $p_0(x)$ does not satisfy the boundary condition
(\ref{boucon1}), it is necessary to construct a boundary layer correction to $p_0(x)$ near $x=A$.

We introduce a stretched variable
\ce
\xi=\frac{x-A}{\e^\beta}
\de
with $\beta>0$ determined later. Defining $F(\xi)=p_0(A+\xi\e^\beta)$ and inserting it into Eq.(\ref{intdif2}), we obtain
\be
b(A+\xi\e^\beta)F'(\xi)\e^{-\beta}+\e^{\alpha-\alpha\beta}\int_{\mR\setminus\{0\}}\big(F(\xi+u)-F(\xi)
-F'(\xi)u\big)\nu(\dif u)=0.
\label{intdifA}
\ee
Set $-\beta=\alpha-\alpha\beta$. That is, we take $\beta=\frac{\alpha}{\alpha-1}$. Multiplying
Eq.(\ref{intdifA}) with $\e^{\beta}$ and letting $\e\rightarrow0$, we get
\be
b(A)F'(\xi)+\int_{\mR\setminus\{0\}}\big(F(\xi+u)-F(\xi)
-F'(\xi)u\big)\nu(\dif u)=0,
\label{intdifAc}
\ee
with the boundary condition
\be
F(\xi)=0, \quad \xi\leq0,
\label{matcon1}
\ee
and the matching condition
\be
\lim\limits_{\xi\rightarrow\infty}F(\xi)=1.
\label{matcon2}
\ee

By \cite{ja2}, we know that  the system (\ref{intdifAc})-(\ref{matcon2}) is
solvable, although the solution cannot be expressed in terms of elementary functions. For the special example
we consider in Section \ref{example}, the boundary layer function will be constructed explicitly. So,
$$
p_0(x)=F\left(\frac{x-A}{\e^{\frac{\alpha}{\alpha-1}}}\right).
$$

Thus an asymptotic solution of $p(x)$ is, for sufficiently small $\e$,
\ce
p(x)\approx F\left(\frac{x-A}{\e^{\frac{\alpha}{\alpha-1}}}\right).
\de

{\bf Case 2:}   $b(x)<0$ for $x\in(A,B)$. Again, in this case the deterministic dynamical system $\dot{x}=b(x)$ has no equilibrium states and all orbits move to the left.

Thus as $\e\rightarrow0$, $p(x)\rightarrow 0$.
We assume that $p(x)$ has the following   expansion
\be
p(x)=p_0(x)+\e^\alpha p_1(x)+\e^{2\alpha}p_2(x)+\cdots.
\label{asym3}
\ee
Similar to {\bf Case 1}, we obtain the leading-order equation for $p_0(x)$
\ce
b(x)p_0^{'}(x)=0, \qquad x\in(A,B).
\de
So, $p_0(x)=0$ for $x\in(A,B)$, because $p(x)\rightarrow 0$ as $\e\rightarrow0$. Since $p_0(x)$ does not satisfy the boundary condition
(\ref{boucon2}), it is necessary to construct a boundary layer correction to $p_0(x)$ near $x=B$.

We introduce a stretched variable
\ce
\varsigma=\frac{B-x}{\e^\beta},
\de
where $\beta$ is the same as one in {\bf Case 1}. Defining $G(\varsigma)=p_0(B-\varsigma\e^\beta)$ and inserting it into Eq.(\ref{intdif2}), we obtain
\be
 -b(B-\varsigma\e^\beta)G'(\varsigma)\e^{-\beta}+\e^{\alpha-\alpha\beta}\int_{\mR\setminus\{0\}}\big[G(\varsigma-u)-G(\varsigma)
 -G'(\varsigma)(-u)\big]\nu(\dif u)=0.
\label{intdifB}
\ee
Multiplying Eq.(\ref{intdifB}) with $\e^{\beta}$ and letting $\e\rightarrow0$,   we obtain
\be
-b(B)G'(\varsigma)+\int_{\mR\setminus\{0\}}\big(G(\varsigma+u)-G(\varsigma)
-G'(\varsigma)u\big)\nu(\dif u)=0,
\label{intdifBc}
\ee
with the boundary condition
\be
G(\varsigma)=1, \quad \varsigma\leq0,
\label{matcon3}
\ee
and the matching condition
\be
\lim\limits_{\varsigma\rightarrow\infty}G(\varsigma)=0.
\label{matcon4}
\ee

By \cite{ja2}, the system (\ref{intdifBc})-(\ref{matcon4})
is solvable, although the solution cannot be expressed in terms of elementary functions. So,
$$
p_0(x)=G\left(\frac{B-x}{\e^{\frac{\alpha}{\alpha-1}}}\right).
$$

Thus we obtain   an asymptotic solution of Eq.(\ref{intdif2})
\ce
p(x)\approx G\left(\frac{B-x}{\e^{\frac{\alpha}{\alpha-1}}}\right).
\de

{\bf Case 3:} There exists only one $\bar{x}\in(A,B)$ such that $b(\bar{x})=0$ and $b'(\bar{x})>0$
 (Assume that  $b$ is differentiable at $\bar{x}$).
In this case the deterministic dynamical system $\dot{x}=b(x)$ has one \emph{unstable} equilibrium state $\bar{x}$.
Then as $\e\rightarrow0$, $p(x)\rightarrow1$ for $\bar{x}<x\leq B$ and
$p(x)\rightarrow0$ for $A\leq x<\bar{x}$. We assume that $p(x)$ has the following expansion
\be
p(x)=p_0(x)+\e^\alpha p_1(x)+\e^{2\alpha}p_2(x)+\cdots.
\label{asym4}
\ee
As in {\bf Case 1}, we obtain that the leading-order
equation for $p_0(x)$ is
\ce
b(x)p_0^{'}(x)=0, \qquad x\in(A,B).
\de
So,
$$
p_0(x)=\left\{\begin{array}{l}
1, \quad \bar{x}<x\leq B,\\
0, \quad A\leq x<\bar{x}.
\end{array}
\right.
$$
Although $p_0(x)$ partially satisfies the `exterior' conditions (\ref{boucon1}) and (\ref{boucon2}),
the value of $p_0(x)$ around $\bar{x}$ is unknown. Therefore, it is necessary to construct
an internal boundary layer correction to $p_0(x)$ near $x=\bar{x}$.

We introduce a stretched variable
\ce
\eta=\frac{x-\bar{x}}{\e}.
\de
Define $H(\eta)=p_0(\bar{x}+\eta\e)$ and insert it into Eq.(\ref{intdif2}).
Then Eq.(\ref{intdif2}) becomes
\ce
b(\bar{x}+\eta\e)H'(\eta)\e^{-1}+\int_{\mR\setminus\{0\}}\big(H(\eta+u)-H(\eta)
-H'(\eta)u\big)\nu(\dif u)=0.
\de
Letting $\e\rightarrow0$ and using the L'Hospital's rule, we get
\be
b'(\bar{x})\eta H'(\eta)+\int_{\mR\setminus\{0\}}\big(H(\eta+u)-H(\eta)
-H'(\eta)u\big)\nu(\dif u)=0,
\label{intdifCc}
\ee
with the matching conditions
\be
\lim\limits_{\eta\rightarrow-\infty}H(\eta)=0, \label{matcon5}\\
\lim\limits_{\eta\rightarrow\infty}H(\eta)=1.
\label{matcon6}
\ee

By \cite{ja2}, Eq.(\ref{intdifCc}) is solvable. So,
$$
p_0(x)=H\left(\frac{x-\bar{x}}{\e}\right).
$$

Thus we obtain an asymptotic solution of Eq.(\ref{intdif2})
\ce
p(x)\approx H\left(\frac{x-\bar{x}}{\e}\right).
\de

{\bf Case 4:} There exists only one $\bar{x}\in(A,B)$ such that $b(\bar{x})=0$ and $b'(\bar{x})<0$
(Assume that  $b$ is differentiable at $\bar{x}$).
In this case the deterministic dynamical system $\dot{x}=b(x)$
has one \emph{stable} equilibrium state $\bar{x}$.
We assume that $p(x)$ has the following   expansion
\be
p(x)=p_0(x)+\e^\alpha p_1(x)+\e^{2\alpha}p_2(x)+\cdots.
\label{asym5}
\ee
As in {\bf Case 1}, we obtain   the leading-order
equation for $p_0(x)$
\ce
b(x)p_0^{'}(x)=0, \qquad x\in(A,B).
\de
So, $p_0(x)=C$ for $x\in(A,B)$. Because of not knowing at which endpoint there
will be a boundary layer correction, we   construct asymptotic approximations
near both endpoints. If there is a boundary layer correction near $x=A$ and $x=B$,
respectively,  as in {\bf Case 1} and {\bf Case 2},
we get near $x=A$
$$
p_0(x)=C\cdot F\left(\frac{x-A}{\e^{\frac{\alpha}{\alpha-1}}}\right),
$$
and near $x=B$
$$
p_0(x)=C+\left(1-C\right)\cdot G\left(\frac{B-x}{\e^{\frac{\alpha}{\alpha-1}}}\right).
$$
Thus we have an asymptotic solution
$$
p(x)\approx C\cdot F\left(\frac{x-A}{\e^{\frac{\alpha}{\alpha-1}}}\right)
+\left(1-C\right)\cdot G\left(\frac{B-x}{\e^{\frac{\alpha}{\alpha-1}}}\right).
$$
Since $F, G$ cannot be expressed in terms of elementary functions,
it is hard to determine $C$. But under the condition that $b(A)b(B)\neq0$,
i.e., $A$ and $B$ are not equilibrium states for $\dot{x}=b(x)$, for a concrete example
in the next section, Example \ref{example5.3}, we introduce a method to determine  the value of $C$.

The following graphs present the movement directions of solution orbits for the deterministic
dynamical system $\dot{x}=b(x)$ in above four cases:

\begin{figure}
{\scalebox{0.8}{\includegraphics*[-8,-2][1070,180]{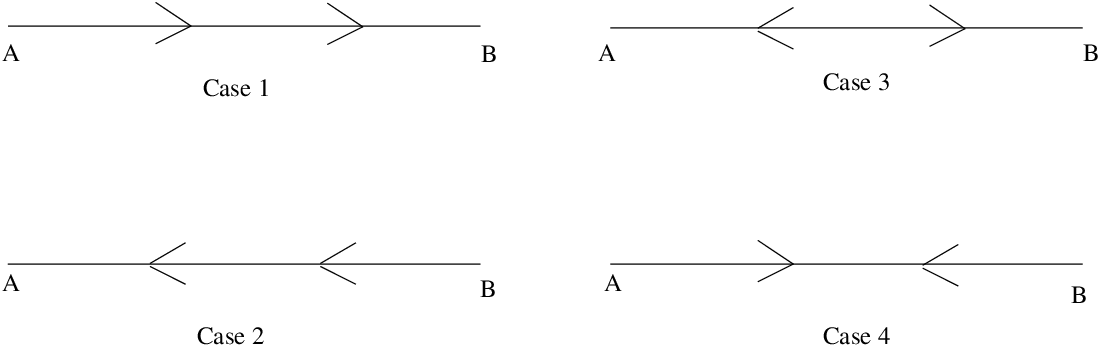}}}\\

\vspace*{5mm}

\caption{The movement directions of solution orbits for the
deterministic dynamical system $\dot{x}=b(x)$ in Case 1, 2, 3, 4.}
\label{fig1}
\end{figure}

\br\label{othercases} For other cases with more than one equilibrium
state for the deterministic dynamical system $\dot{x}=b(x)$, they
are so complex that, if we assume that $p(x)$ has the following
expansion \ce p(x)=p_0(x)+\e^\alpha
p_1(x)+\e^{2\alpha}p_2(x)+\cdots, \de boundary layer analysis of
$p_0(x)$ couldn't be done and then the asymptotic solution for
Eq.(\ref{intdif2}) isn't given explicitly. Therefore, we don't
consider these cases. \er

\section{Examples}
\label{example}

In this section we consider three examples.
Example \ref{example5.1}, Example \ref{example5.2} and Example \ref{example5.3} correspond to  our methods in Section
\ref{sde-contjump}, {\bf Case 1} and {\bf Case 2} of Section \ref{sde-jump}, respectively.

\begin{example} \label{example5.1}
Consider the following scalar \textsc{SDE} with a Brownian motion and a symmetric
$\alpha$-stable L\'evy motion:
 \ce\left\{\begin{array}{l}
\dif X_{t}=\dif W_t+\e\dif L_t^{\alpha},\\
X_0=x.
\end{array}
\right.
\de
The unique solution is denoted as $X_t(x)$. We take $(A,B)=(-1,1)$ and $[B,\infty)=[1,\infty)$.
The escape probability $p(x)$, for $X_{t}(x)$ with $0<\alpha<1$ from $(-1,1)$ to $[1,\infty)$,
satisfies the following integro-differential equation
\be
\frac{1}{2}p^{''}(x)+\e^\alpha\int_{\mR\setminus\{0\}}\left(p(x+u)-p(x)\right)\nu(\dif u)=0, \quad x\in(-1,1),
\label{side}
\ee
and the exterior conditions
\ce
p(x)|_{(-\infty,-1]}=0, \qquad p(x)|_{[1,\infty)}=1.
\de
We seek an asymptotic solution of $p(x)$ as follows
\ce
p(x)\approx p_0(x)+\e^\alpha p_1(x),
\de
where
\ce
p_0(x)=\left\{\begin{array}{l}
0, \qquad x\leq -1,\\
\frac{x+1}{2}, \quad -1<x<1,\\
1, \qquad x\geq 1,
\end{array}
\right.
\de
and
\ce
p_1(x)=\left\{\begin{array}{l}
0, \qquad\qquad\qquad\qquad\qquad\qquad\quad\qquad\qquad\qquad\qquad\qquad\qquad\quad x\leq -1,\\
\frac{C_{1,\alpha}}{(-\a)(1-\a)(2-\a)(3-\a)}\big[(1-x)^{3-\a}-2^{3-\a}+(3-\a)2^{2-\a}(x+1)-(1+x)^{3-\a}\big]\\
-\frac{x+1}{2}\frac{C_{1,\alpha}}{(-\a)(2-\a)(3-\a)}2^{3-\a}+\frac{x+1}{2}, \qquad\qquad\qquad\qquad\qquad\qquad\quad -1<x<1,\\
1, \qquad\qquad\qquad\qquad\qquad\qquad\quad\qquad\qquad\qquad\qquad\qquad\qquad\qquad x\geq 1.
\end{array}
\right.
\de

By the same calculation as above, we could obtain the asymptotic solution of $p(x)$ with $1<\alpha<2$.

Next we use the numerical method in \cite{GaoDuan} to study Eq.(\ref{side}).

\begin{figure}[H]
\begin{center}
\includegraphics{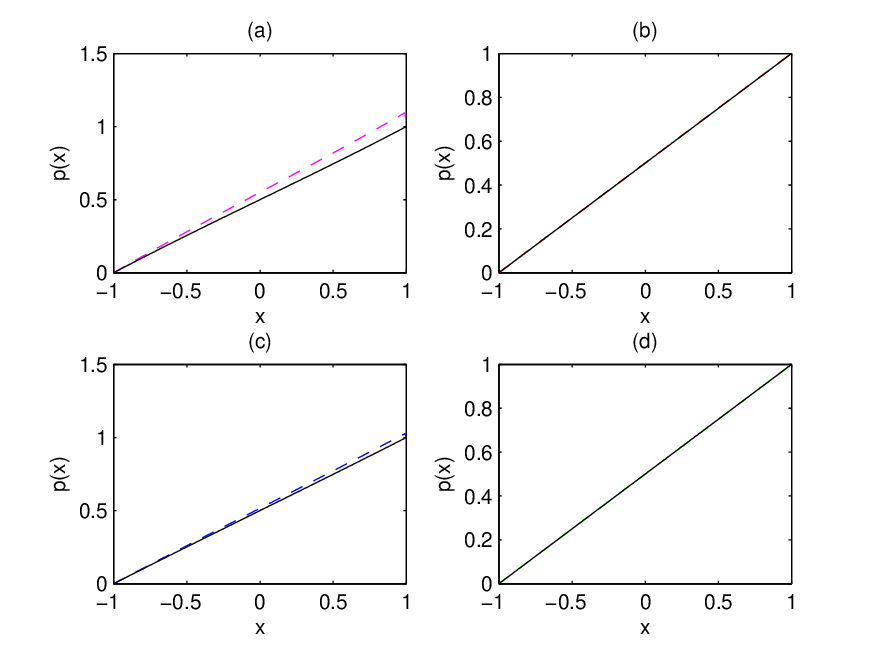}
\end{center}
\caption{Comparison between the asymptotic solution and the
numerical solution of Eq.(\ref{side}) for small $\e$. (a) $\a=0.5,
\e=0.01$. The asymptotic solution is shown with dashed line while
the numerical solution is displayed with solid line. (b) $\a=1.5,
\e=0.01$. (c) $\a=0.5, \e=0.001$. (d) $\a=1.5, \e=0.001$.}
\label{fig2}
\end{figure}

Figure \ref{fig2} shows that if $\e$ is fixed and $\a$ turns large, the difference
between the asymptotic solution and the numerical solution of Eq.(\ref{side}) will
become small; if $\a$ is fixed and $\e$ becomes large, the difference will turn large, because the asymptotic solution is for sufficiently small $\e$.
\end{example}

\begin{example} \label{example5.2}
For the deterministic dynamical system
\ce
\dot{x}=x(1-\theta x)-\beta\frac{x}{x+1},
\de
where $0<\theta<1$, $1<\beta<\frac{(\theta+1)^2}{4\theta}$ and the potential function is
\ce
U(x)=-\frac{x^2}{2}+\frac{\theta x^3}{3}+\beta x-\beta\ln(x+1),
\de
$x$ is the normalized molecular density of tumor cells with respect to the maximum tissue
capacity (\cite{fosg}). The system has two stable states and one unstable state:
\ce
&&x_1=0,\\
&&x_2=\frac{1-\theta-\sqrt{(1-\theta)^2-4\theta(\beta-1)}}{2\theta},\\
&&x_3=\frac{1-\theta+\sqrt{(1-\theta)^2-4\theta(\beta-1)}}{2\theta}.
\de

Without random fluctuations, system states finally approach one of the two
stable states: (i) either the stable state $x_1 = 0$, where no tumor cells are
present, namely, the tumor-free state (or the state of tumor extinction), (ii)
or the other stable state $x_3$, where the tumor cell density does not increase
but stays at a certain constant level, namely, the state of stable tumor.

Under the environment fluctuations, the tumor density is described by the following
scalar \textsc{SDE}  with a symmetric $\alpha$-stable L\'evy motion, with $1<\a<2$,
\ce\left\{\begin{array}{l}
\dif X_{t}=\left[X_t(1-\theta X_t)-\beta\frac{X_t}{X_t+1}\right]\,\dif t+\e\dif L_t^{\alpha},\\
X_0=x.
\end{array}
\right.
\de
The unique solution is denoted as $X_t(x)$. We take $(A,B)=(x_1,x_3)$ and $(-\infty,A]=(-\infty,x_1]$.
The escape probability $p(x)$, for $X_{t}(x)$ from $(x_1,x_3)$ to $(-\infty,x_1]$, i.e. the likelihood
of tumor extinction, satisfies
the following integro-differential equation
\ce
\left[x(1-\theta x)-\beta\frac{x}{x+1}\right]p'(x)+\e^\alpha\int_{\mR\setminus\{0\}}\left(p(x+u)-p(x)
-p'(x)u\right)\nu(\dif u)=0,\\
x\in(x_1,x_3),
\de
with the exterior conditions
\ce
p(x)|_{(-\infty,x_1]}=1, \qquad p(x)|_{[x_3,\infty)}=0.
\de
Since $b(x)=x(1-\theta x)-\beta\frac{x}{x+1}, b(x_2)=0, b'(x_2)>0$, by the result of {\bf Case 3} in Section \ref{sde-jump},
an asymptotic solution of $p(x)$ is given by
\ce
p(x)\approx \bar{H}\left(\frac{x-x_2}{\e}\right),
\de
where $\bar{H}(x)$ solves the following equation
\ce
b'(x_2)\eta \bar{H}'(\eta)+\int_{\mR\setminus\{0\}}\big(\bar{H}(\eta+u)-\bar{H}(\eta)
-\bar{H}'(\eta)u\big)\nu(\dif u)=0,
\de
with the matching conditions
\ce
\lim\limits_{\eta\rightarrow-\infty}\bar{H}(\eta)=1, \\
\lim\limits_{\eta\rightarrow\infty}\bar{H}(\eta)=0.
\de
\end{example}

\begin{example} \label{example5.3}
Consider the following scalar \textsc{SDE}  with a symmetric
$\alpha$-stable L\'evy motion, with $1<\a<2$,
\ce\left\{\begin{array}{l}
\dif X_{t}=-X_t\,\dif t+\e\dif L_t^{\alpha},\\
X_0=x.
\end{array}
\right.
\de
The unique solution is denoted as $X_t(x)$. We take $(A,B)=(-1,1)$ and $[B,\infty)=[1,\infty)$.
The escape probability $p(x)$, for $X_{t}(x)$ from $(-1,1)$ to $[1,\infty)$, satisfies
the following integro-differential equation
\be
-x p'(x)+\e^\alpha\int_{\mR\setminus\{0\}}\left(p(x+u)-p(x)
-p'(x)u\right)\nu(\dif u)=0,\no\\
x\in(-1,1),
\label{spid}
\ee
with the exterior conditions
\be
p(x)|_{(-\infty,-1]}=0, \qquad p(x)|_{[1,\infty)}=1.
\label{spbc}
\ee
Since $b(x)=-x, b(0)=0, b'(0)<0, b(-1)b(1)=-1\neq0$, by the result of {\bf Case 4} in Section \ref{sde-jump},
an asymptotic solution of $p(x)$ is given by
\be
p(x)\approx C\cdot F\left(\frac{x+1}{\e^{\frac{\alpha}{\alpha-1}}}\right)
+(1-C)\cdot G\left(\frac{1-x}{\e^{\frac{\alpha}{\alpha-1}}}\right).
\label{asymso}
\ee

Specially, take the L\'evy measure
$$
\nu(\dif u)=\frac{\kappa}{|u|^{1+\alpha}}\cdot 1_{|u|\leq1}\dif u,
$$
where $\kappa>0$ is a constant(\cite{Woy}). Thus, the function $F$ can be given explicitly
by
$$
F(x)=\left\{\begin{array}{l}
1-e^{-\gamma x}, \quad x>0,\\
0, \qquad\qquad x\leq0,
\end{array}
\right.
$$
where $\gamma>0$ satisfies the following integral equation
$$
\gamma-\int_{-1}^1\left(e^{-\gamma u}-1-(-\gamma)u\right)\frac{\kappa}
{|u|^{1+\alpha}}\dif u=0.
$$
By the relation between $F$ and $G$, we can obtain
$$
G(x)=\left\{\begin{array}{l}
e^{-\gamma x}, \qquad x>0,\\
1, \qquad\quad x\leq0.
\end{array}
\right.
$$
So, the asymptotic solution of Eq.(\ref{spid}) is given by
$$
p(x)\approx C\left(1-\exp\left\{-\gamma\left(\frac{x+2}{\e^{\frac{\alpha}{\alpha-1}}}\right)\right\}\right)
+(1-C)\exp\left\{-\gamma\left(\frac{1-x}{\e^{\frac{\alpha}{\alpha-1}}}\right)\right\}.
$$

To determine $C$, we multiply Eq.(\ref{spid}) by the solution $\rho(x)$ of the steady  Fokker-Planck equation
\be
-(-x\rho(x))'+\e^\alpha\int_{\mR\setminus\{0\}}\left(\rho(x+u)-\rho(x)
-\rho'(x)u\right)\nu(\dif u)=0
\label{fpe}
\ee
and integrate over $(-1,1)$, to obtain
\be
\int_{-1}^1\big(-x\rho(x)\big)p'(x)\dif x+\e^\alpha\int_{-1}^1\rho(x)\dif x\int_{\mR\setminus\{0\}}\left(p(x+u)-p(x)
-p'(x)u\right)\nu(\dif u)=0.
\label{info}
\ee
To (\ref{info}), by integration by parts and using (\ref{fpe}), we get
\ce
&&-\rho(1)-\e^\alpha\int_{-1}^1p(x)\dif x\int_{\mR\setminus\{0\}}\left(\rho(x+u)-\rho(x)
-\rho'(x)u\right)\nu(\dif u)\\
&&+\e^\alpha\int_{-1}^1\rho(x)\dif x\int_{\mR\setminus\{0\}}\left(p(x+u)-p(x)
-p'(x)u\right)\nu(\dif u)=0.
\de

Applying Cauchy principal value and (\ref{spbc}), we have
\be
&&-\frac{\e^{-\alpha}\rho(1)}{C_{1,\a}}+\int_{-1}^1\rho(x)\frac{(1-x)^{-\alpha}}{\alpha}\dif x\no\\
&=&\int_{-1}^1p(x)\dif x\left[\int_{-\infty}^{-1}\frac{\rho(y)}{|y-x|^{1+\a}}\dif y
+\int_1^{\infty}\frac{\rho(y)}{|y-x|^{1+\a}}\dif y\right].
\label{capr}
\ee
By \cite[Proposition 3.2]{arw}, the Fourier transform of $\rho(k)$ is given by
\ce
\hat{\rho}(k)=\exp\{-\frac{\e^\alpha}{\a}|k|^\alpha\}.
\de
Replacing $\rho(x)$ and $p(x)$ by $\frac{1}{2\pi}\int_{\mR}e^{ixk}\hat{\rho}(k)\dif k$
and (\ref{asymso}), respectively, and letting $\e\rightarrow0$, we can obtain $C$ from (\ref{capr}).
\end{example}

\bigskip

\emph{Acknowledgement}. We thank Xingye Kan for help with Example \ref{example5.1}.
This work was done while Huijie Qiao was visiting the Institute for Pure and Applied Mathematics (IPAM), Los Angeles.

\end{document}